\documentclass[11pt]{article}
\usepackage{amsmath}
\usepackage{amssymb}

\font\tenmsb=msbm10

\def\eps{\varepsilon}
\font\tencmmib=cmmib10 \skewchar\tencmmib '60
\newfam\cmmibfam
\textfont\cmmibfam=\tencmmib

\font\tenmsb=msbm10

\def\Bbb#1{\hbox{\tenmsb#1}}

\def\bbox{\quad\hbox{\vrule \vbox{\hrule \vskip2pt \hbox{\hskip2pt
\vbox{\hsize=1pt}\hskip2pt} \vskip2pt\hrule}\vrule}}
\def\lessim{\ \lower4pt\hbox{$
\buildrel{\displaystyle <}\over\sim$}\ }
\def\gessim{\ \lower4pt\hbox{$\buildrel{\displaystyle >}
\over\sim$}\ }

\def\eps{{\varepsilon}}

\def\Bbb E{\mathbb{E}}
\def\Bbb R{\mathbb{R}} 
\newtheorem{lemma}{Lemma}
\newtheorem{theorem}{Theorem}
\newtheorem{corollary}{Corollary}
\makeatletter
\@addtoreset{equation}{section}

\makeatother

\font\tencmmib=cmmib10 \skewchar\tencmmib '60
\newfam\cmmibfam
\textfont\cmmibfam=\tencmmib

\font\tenmsb=msbm10

\def\Bbb#1{\hbox{\tenmsb#1}}

\def\bbox{\quad\hbox{\vrule \vbox{\hrule \vskip2pt \hbox{\hskip2pt
\vbox{\hsize=1pt}\hskip2pt} \vskip2pt\hrule}\vrule}}
\def\lessim{\ \lower4pt\hbox{$
\buildrel{\displaystyle <}\over\sim$}\ }
\def\gessim{\ \lower4pt\hbox{$\buildrel{\displaystyle >}
\over\sim$}\ }

\def\EE{{\Bbb E} }

\def\eps{\varepsilon}

\def\go0{\to 0}

\def\leftitem#1{\item{\hbox to\parindent{\enspace#1\hfill}}}

\def\qed{{$\hfill \bbox$}}

\def\sg{\sigma}

\def\sg2{\sigma^2}

\def\__{_{\infty}}
\begin{document}
{\baselineskip=16.5pt

\title{\bf Sharp Oracle Inequalities in Low Rank Estimation}

\author{
{\bf Vladimir Koltchinskii}
\thanks{Partially supported by NSF Grants DMS-1207808, DMS-0906880 and CCF-0808863}
\\ School of Mathematics
\\ Georgia Institute of Technology
\\ Atlanta, GA 30332-0160
\\ vlad@math.gatech.edu
}

\maketitle

\begin{abstract}
The paper deals with the problem of penalized empirical risk minimization 
over a convex set of linear functionals 
on the space of Hermitian matrices with convex loss and nuclear norm 
penalty. Such penalization is often used in low rank matrix recovery 
in the cases when the target function 
can be well approximated by a linear functional generated by a 
Hermitian matrix of relatively small rank (comparing with the size 
of the matrix).
Our goal is to prove sharp low rank oracle inequalities that involve 
the excess risk (the approximation error) with constant equal to one
and the random error term with correct 
dependence on the rank of the oracle. 
\end{abstract}

%{{\bf Keywords and phrases:}
%low rank matrix estimation, 
%von Neumann entropy, 
%matrix regression, 
%empirical processes, 
%noncommutative Bernstein inequality, 
%quantum state tomography}

%{{\bf 2010 AMS Subject Classification:} 62J99, 62H12, 60B20, 60G15, 81Q99}

\medskip

%\section{Introduction}

%\medskip

\section{Main Result}

Let $(X,Y)$ be a couple, where $X$ is a random variable in the space ${\mathbb H}_m$ of $m\times m$ Hermitian matrices and $Y$ is a random response variable with values in a Borel subset $T\subset {\mathbb R}.$ Let $P$ be the distribution of $(X,Y)$ and 
let $\Pi$ denote the marginal distribution of $X.$
The goal is to predict $Y$ based on an observation of $X.$ More precisely, let $\ell : T\times {\mathbb R}\mapsto {\mathbb R}_+$ be 
a measurable loss function. We will assume in what follows that, for all $y\in T,$ $\ell(y;\cdot)$ is convex. Given a measurable function $f:{\mathbb H}_m\mapsto {\mathbb R}$ (a ``prediction rule''), denote $(\ell \bullet f)(x,y):=\ell (y;f(x))$ and define the risk of $f$ as 
$$
P(\ell \bullet f)={\mathbb E}\ell (Y;f(X)).
$$
Then, one can view the prediction problem as risk minimization: the goal is to find 
a function $f_{\ast}:{\mathbb H}_m\mapsto {\mathbb R}$ that minimizes the risk 
$P(\ell \bullet f)$ over the class of all measurable prediction rules $f:{\mathbb H}_m\mapsto {\mathbb R}$ (provided that such a function exists), or, more realistically, to find a reasonably good approximation of $f_{\ast}.$ To this 
end, one wants to find a function $f$ for which the excess risk 
${\cal E}(f):=P(\ell \bullet f)-\inf_{g:{\mathbb H}_m\mapsto {\mathbb R}}P(\ell \bullet g)$ is small enough. Of course, the risk $P(\ell\bullet f)$ depends on 
the distribution $P$ of $(X,Y),$ which is, most often, unknown. In such cases,
the problem has to be solved based on the training data $(X_1,Y_1),\dots ,(X_n,Y_n)$
that consists of $n$ independent copies of $(X,Y).$ We will be especially 
interested in the problems in which matrices are large and the optimal 
prediction rule $f_{\ast}$ can be well approximated by a linear function $f_{S}(\cdot):=\langle S,\cdot \rangle,$ where $S\in {\mathbb H}_m$ is a low 
rank Hermitian matrix, that is, when there exists a low rank matrix $S$ (an oracle)
such that the excess risk ${\cal E}(f_S)$ is small. Here and in what follows, $\langle \cdot,\cdot \rangle$ denotes the Hilbert-Schmidt (Frobenius) inner 
product in ${\mathbb H}_m.$
In such problems, we would like 
to find an estimator $\hat S$ based on the training data $(X_1,Y_1),\dots, (X_n,Y_n)$
such that the excess risk ${\cal E}(f_{\hat S})$ of the estimator can be bounded from above 
by the excess risk ${\cal E}(f_{S})$ of an arbitrary oracle $S\in {\mathbb H}_m$
plus an error term that properly depends on the rank of the oracle. The resulting 
bounds on the excess risk ${\cal E}(f_{\hat S})$ of the estimator $\hat S$ are  supposed to hold with a guaranteed high probability and they are often called 
``low rank oracle inequalities.'' We will consider below rather traditional 
estimator $\hat S$ based on penalized empirical risk minimization with a nuclear 
norm penalty:
\begin{equation}
\label{nuclear_ERM}
\hat S:= {\rm argmin}_{S\in {\mathbb D}}\Bigl[P_n(\ell \bullet f_{S})+\eps \|S\|_1\Bigr],
\end{equation}
where ${\mathbb D}\subset {\mathbb H}_m$ is a closed convex set, $0\in {\mathbb D},$ $P_n$ is the 
empirical distribution based on the training data $(X_1,Y_1),\dots, (X_n,Y_n)$
and 
$$
P_n(\ell \bullet f_S)=n^{-1}\sum_{j=1}^n \ell(Y_j;f_S(X_j))
$$ 
is the corresponding empirical risk with respect to the loss $\ell,$ 
$\|S\|_1:={\rm tr}(|S|)={\rm tr}(\sqrt{S^2})$ is the nuclear norm of 
$S$ and $\eps\geq 0$ is the regularization parameter. Clearly, optimization 
problem (\ref{nuclear_ERM}) is convex. In fact, it is a standard convex 
relaxation of penalized empirical risk minimization with a penalty 
proportional to the rank of $S,$ denoted in what follows by ${\rm rank}(S),$ which would not be a computationally tractable problem. Such convex relaxations have been extensively studied in the recent years 
(see Recht, Fazel and Parrilo (2010), Candes and Recht (2009), Candes and Tao (2010), Candes and Plan (2011),
Gross (2011), Rohde and Tsybakov (2011), Negahban and Wainwright (2010), Koltchinskii (2011), Koltchinskii, Lounici and Tsybakov (2011) 
and references therein).

To state our main result (a sharp low rank oracle inequality for the estimator $\hat S$), we first introduce some assumptions and notations. In what follows, assume 
that for some constant $a>0,$ $|\langle S,X\rangle|\leq a\ {\rm a.s.},\ S\in {\mathbb D}.$ It will be also assumed that $\ell$ is a convex loss of \it quadratic 
type. \rm More precisely, suppose that, for all $y\in T,$ $\ell (y,\cdot)$ is 
twice continuously differentiable convex function in $[-a,a]$ with 
$Q:=\sup_{y\in T}\ell(y;0)<+\infty,$
$$
L(a):= \sup_{y\in T}\sup_{u\in [-a,a]}
\Bigl[|\ell'(y;0)|+\ell^{''}(y;u)a\Bigr]<+\infty
$$
and
$$
\tau(a):= \inf_{y\in T}\inf_{u\in [-a,a]} \ell^{''}(y;u)>0.
$$
Here $\ell^{'}, \ell^{''}$ denote the first and the second derivatives of 
the loss $\ell(y,u)$ with respect to $u.$ Many important losses in regression 
and in large margin classification problems are of quadratic type. In particular,
if $\ell(y;u)=(y-u)^2, y,u\in [-a,a]$ (regression with quadratic loss and with 
bounded response), then $L(a)=4a$ and $\tau(a)=2.$ Exponential loss $\ell(y,u)=
e^{-yu}, y\in \{-1,1\}, u\in [-a,a]$ often used in large margin methods for binary 
classification is also of quadratic type. 

In what follows, $\|\cdot\|_2$ denotes the Hilbert--Schmidt (Frobenius) norm 
of Hermitian matrices (generated by the inner product $\langle \cdot,\cdot \rangle$)
and $\|\cdot\|$ denotes the operator norm.

We will use certain characteristics of matrices $S\in {\mathbb D}$ 
that are related to matrix versions of restricted isometry property
(see, e.g., Koltchinskii (2011), Chapter 9 and references therein).
Let $S\in {\mathbb D}$ be a matrix with spectral representation 
$S=\sum_{j=1}^r \lambda_j (\phi_j\otimes \phi_j),$ where $r:={\rm rank}(S),$ $\lambda_j$ are 
non-zero eigenvalues of $S$ (repeated with their multiplicities) and 
$\phi_j\in {\mathbb C}^m$ are the corresponding orthonormal eigenvectors. In what follows,
we denote 
$$
{\rm sign}(S):=\sum_{j=1}^r {\rm sign}(\lambda_j)(\phi_j\otimes \phi_j),\ \
L:={\rm supp}(S):={\rm l.s.}(\phi_1,\dots, \phi_r).
$$
Let ${\cal P}_L, {\cal P}_L^{\perp}$ be the following orthogonal 
projectors in the space $({\mathbb H}_m,\langle \cdot,\cdot\rangle):$
$$
{\cal P}_L (A):= A-P_{L^{\perp}}AP_{L^{\perp}},\  
{\cal P}_L^{\perp} (A):= P_{L^{\perp}}AP_{L^{\perp}},\ A\in {\mathbb H}_m 
$$
(here $L^{\perp}$ is the orthogonal complement of $L$). Clearly, 
we have $A={\cal P}_L A+ {\cal P}_L^{\perp}A, A\in {\mathbb H}_m,$
providing a decomposition of a matrix $A$ into a ``low rank part''
${\cal P}_L A$ and a ``high rank part'' ${\cal P}_L^{\perp} A.$
Given $b>0,$ define the following cone in the space ${\mathbb H}_m$
$$
{\cal K}({\mathbb D};L;b):= 
\Bigl\{A\in {\rm l.s.}({\mathbb D}): \|{\cal P}_L^{\perp}(A)\|_1\leq b\|{\cal P}_L (A)\|_1\Bigr\}
$$
that consists of matrices $A$ with a ``dominant'' low rank part. 
Let 
$$
\beta^{(b)}({\mathbb D};L;\Pi):=
\inf\Bigl\{\beta>0: \|{\cal P}_L (A)\|_2 \leq \beta \|f_A\|_{L_2(\Pi)}, A\in {\cal K}({\mathbb D};L;b)\Bigr\}.
$$
This quantity is known to be bounded from above by a constant in the case 
when the matrix form of ``distribution dependent'' restricted isometry condition holds for $r=4{\rm rank}(S)$
(see Koltchinskii (2011), Section 9.1). 
In what follows, we will use the following characteristic of oracle $S:$
$$
\beta(S):= \beta^{(5)}({\mathbb D};L;\Pi),\ L:={\rm supp}(S).
$$

For arbitrary $t>0$ and $S\in {\mathbb D},$ denote
$$
t(S;\eps):=
t+3\log\Bigl(B\log_2 \Bigl(\|S\|_1\vee n \vee \eps\vee Q\vee a^{-1}\vee (L(a))^{-1}\vee 2\Bigr)\Bigr),
$$
where $B>0$ is a constant.
Let 
$$
\Delta:= {\mathbb E}\biggl\|\frac{1}{\sqrt{n}}\sum_{j=1}^n \eps_j X_j\biggr\|,
$$
where $\{\eps_j\}$ are i.i.d. Rademacher random variables independent 
of $\{X_j\}.$

\begin{theorem}
\label{main_theorem}
There exist a numerical constant $B>0$ in the definition of $t(S;\eps)$ and 
numerical constants $C,D>0$ such that for all $t>0$ and all 
\begin{equation}
\label{tr_eps}
\eps\geq \frac{D L(a)\Delta}{\sqrt{n}} ,
\end{equation}
with probability at least $1-e^{-t},$
\begin{equation}
\label{main_bound}
{\cal E}(f_{\hat S})
\leq 
\inf_{S\in {\mathbb D}}\biggl[{\cal E}(f_S)+ 
\Bigl(\frac{3}{\tau(a)}\beta^2(S){\rm rank}(S)\eps^2\bigwedge 2\eps \|S\|_1\Bigr)
+
C(a)\frac{t(S;\eps)}{n}\biggr],
\end{equation}
where
$$
C(a):=C\biggl(\frac{L^2(a)}{\tau(a)}
\bigvee L(a)a\biggr).
$$
\end{theorem}

To control the size of expectation $\Delta$ involved in the 
threshold (\ref{tr_eps}) on $\eps $  one can use a noncommutative 
version of Bernstein inequality due to Ahlswede and Winter (2002). 
Namely, the following upper bound easily follows from this 
inequality (by integrating its exponential tail bounds):
$$
\Delta \leq 4\biggl(\sigma_X\sqrt{\log (2m)}\bigvee U_X \frac{\log(2m)}{\sqrt{n}}\biggr),
$$
where $\sigma_X^2:=\|{\mathbb E} X^2\|$ and 
$U_X:=\Bigl\|\|X\|\Bigr\|_{L_{\infty}}.$ 
This bound can be easily applied to various specific 
sampling models used in low rank matrix recovery, 
such as sampling from an orthonormal basis that 
includes, in particular, matrix completion (see, e.g.,  
Koltchinskii (2011), Chapter 9) leading to more 
concrete results. 

The main feature of oracle inequality (\ref{main_bound}) 
is that it involves the approximation error term ${\cal E}(f_S)$
(the excess risk of the oracle $S$) with constant equal to $1.$
In this sense, bound (\ref{main_bound}) is what is usually called 
a {\it sharp oracle inequality}. Most of low rank oracle inequalities
for the nuclear norm penalization method proved in the recent literature are not sharp in the sense that the oracle excess risk ${\cal E}(f_S)$ is involved in these bounds with a constant strictly larger than $1.$ Sharp oracle inequalities
are especially important in the cases when for all oracles in $S\in {\mathbb D}$ the approximation error is not particularly small.
The first sharp oracle inequalities for nuclear 
norm penalization method were proved in Koltchinskii, Lounici and Tsybakov (2011).
It was done for a ``linearized version'' of least squares method with
nuclear norm penalty. Under the boundedness 
assumption $|\langle S,X\rangle|\leq a\ {\rm a.s.},\ S\in {\mathbb D}$
for some $a>0$ (the same assumption is used in our paper), 
Klopp (2012) proved error bounds (without approximation error term) for the usual matrix LASSO (that is, nuclear norm penalized least squares method). Earlier, Negahban and Wainwright (2010) studied the same problem under additional 
assumptions on the so called ``spikiness'' of the target matrices.
Koltchinskii and Rangel (2012) stated a sharp oracle inequality 
for the same method in the case of noisy matrix completion problem 
with uniform design (in fact, they deduced this result from more general
oracle bounds for estimators of low rank smooth kernels on graphs).   
In the current paper, we establish sharp oracle inequalities for a 
version of the problem with more general losses of quadratic type 
and for general design distributions. 
Note also that the main part of the random error term of bound (\ref{main_bound})
(that is, the term $\frac{3}{\tau(a)}\beta^2(S){\rm rank}(S)\eps^2\bigwedge 2\eps \|S\|_1$) depends correctly on the rank of the oracle. This follows from the  
minimax lower bounds proved in Koltchinskii, Lounici and Tsybakov (2011)
(in fact, the form of the random error term in (\ref{main_bound}) is the same as 
in that paper).

\section{Proof}

We start with the following condition that is necessary 
for $\hat S$ to be a solution
of convex optimization problem (\ref{nuclear_ERM}): for some $\hat V\in \partial \|\hat S\|_1,$
\begin{equation}
\nonumber
P_n(\ell'\bullet f_{\hat S})(f_{\hat S}-f_{S})+ \eps \langle \hat V,\hat S-S\rangle
\leq 0, S\in {\mathbb D}
\end{equation}
(see, e.g., Aubin and Ekeland (1984), Chap. 2, Corollary 6; see also 
Koltchinskii (2011), pp. 198--199). 
This implies that, for all $S\in {\mathbb D}$ 
\begin{equation}
\label{odin}
P(\ell'\bullet f_{\hat S})(f_{\hat S}-f_{S})+ \eps \langle \hat V,\hat S-S\rangle
\leq (P-P_n)(\ell'\bullet f_{\hat S})(f_{\hat S}-f_{S}).
\end{equation}
Since both $\hat S, S\in {\mathbb D},$ we have $|f_{\hat S}(X)|\leq a, |f_S(X)|\leq a$ a.s., and since $\ell$ is a loss of 
quadratic type, it is easy to check that 
\begin{equation}
\label{dva}
P(\ell'\bullet f_{\hat S})(f_{\hat S}-f_{S})\geq P(\ell\bullet f_{\hat S})-P(\ell\bullet f_S)+ \frac{1}{2}\tau (a)\|f_{\hat S}-f_S\|_{L_2(\Pi)}^2.
\end{equation}
If $P(\ell\bullet f_{\hat S})\leq P(\ell\bullet f_S),$ the oracle inequality 
of the theorem holds trivially. So, we assume in what follows that 
$P(\ell\bullet f_{\hat S})> P(\ell\bullet f_S).$
Inequalities (\ref{odin}) and (\ref{dva}) imply that 
\begin{equation}
\label{tri}
P(\ell\bullet f_{\hat S})+\frac{1}{2}\tau (a)\|f_{\hat S}-f_S\|_{L_2(\Pi)}^2
+ \eps \langle \hat V,\hat S-S\rangle
\leq P(\ell\bullet f_S)+(P-P_n)(\ell'\bullet f_{\hat S})(f_{\hat S}-f_{S}).
\end{equation}
The following characterization of subdifferential of the nuclear norm 
is well known:
$$
\partial \|S\|_1=\{{\rm sign}(S)+{\cal P}_L^{\perp}(M): M\in {\mathbb H}_m,\|M\|\leq 1\},
$$
where $L={\rm supp}(S)$ (see, e.g., Koltchinskii (2011), Appendix A.4).
By the duality between the operator and nuclear norms, there 
exists $M\in {\mathbb H}_m$ with $\|M\|\leq 1$ such that 
$$
\langle {\cal P}_L^{\perp}(M), \hat S-S\rangle =
\langle M, {\cal P}_L^{\perp}(\hat S-S)\rangle =
\|{\cal P}_L^{\perp}(\hat S-S)\|_1=\|{\cal P}_L^{\perp}\hat S\|_1.
$$
Then, by monotonicity of subdifferentials of convex functions, we have,
for $V={\rm sign}(S)+{\cal P}_L^{\perp}(M)\in \partial \|S\|_1,$ that 
$$
\langle {\rm sign}(S), \hat S-S\rangle +\|{\cal P}_L^{\perp}\hat S\|_1
=\langle V,\hat S-S\rangle\leq \langle \hat V, \hat S-S\rangle. 
$$
We now substitute the last bound in (\ref{tri}) to get 
\begin{eqnarray}
\label{chetyre}
&&
\nonumber
P(\ell\bullet f_{\hat S})+\frac{1}{2}\tau (a)\|f_{\hat S}-f_S\|_{L_2(\Pi)}^2
+ \eps \|{\cal P}_L^{\perp}\hat S\|_1
\\
&&
\leq P(\ell\bullet f_S)+ \eps\langle {\rm sign}(S), S-\hat S\rangle+
(P-P_n)(\ell'\bullet f_{\hat S})(f_{\hat S}-f_{S}).
\end{eqnarray}

The main part of the proof is a derivation of an upper bound on the empirical 
process $(P-P_n)(\ell'\bullet f_{\hat S})(f_{\hat S}-f_{S}).$ For a given $S\in {\mathbb D}$ and for $\delta_1, \delta_2\geq 0,$ denote 
$$
{\cal A}(\delta_1,\delta_2):=\{A\in {\mathbb D}: A-S\in {\cal K}({\mathbb D};L;5), 
\|f_A-f_S\|_{L_2(\Pi)}\leq \delta_1, \|{\cal P}_L^{\perp}A\|_1\leq \delta_2\},
$$
$$
\tilde {\cal A}(\delta_1,\delta_2,\delta_3):=\{A\in {\mathbb D}: 
\|f_A-f_S\|_{L_2(\Pi)}\leq \delta_1, \|{\cal P}_L^{\perp}A\|_1\leq \delta_2, 
\|{\cal P}_L(A-S)\|_1\leq \delta_3\},
$$
$$
\check {\cal A}(\delta_1,\delta_4):=\{A\in {\mathbb D}: 
\|f_A-f_S\|_{L_2(\Pi)}\leq \delta_1, \|A-S\|_1\leq \delta_4\},
$$
and 
$$
\alpha_n (\delta_1,\delta_2):= \sup\{|
(P_n-P)(\ell'\bullet f_{A})(f_{A}-f_{S})|: A\in {\cal A}(\delta_1,\delta_2)\},
$$
$$
\tilde \alpha_n (\delta_1,\delta_2,\delta_2):= \sup\{|
(P_n-P)(\ell'\bullet f_{A})(f_{A}-f_{S})|: A\in \tilde {\cal A}(\delta_1,\delta_2,\delta_3)\}.
$$
$$
\check \alpha_n (\delta_1, \delta_4):= \sup\{|
(P_n-P)(\ell'\bullet f_{A})(f_{A}-f_{S})|: A\in \check {\cal A}(\delta_1,\delta_4)\}.
$$

\begin{lemma}
\label{main_lemma}
Suppose $0<\delta_k^{-}<\delta_k^{+}, k=1,2,3,4.$ 
Let $t>0$ and 
\begin{eqnarray}
\nonumber
&&
\bar t := t+\sum_{k=1}^2 \log\Bigl([\log_2(\delta_k^+/\delta_k^-)]+2\Bigr)+\log 3, 
\\
&&
\nonumber
\tilde t := t+\sum_{k=1}^3 \log\Bigl([\log_2(\delta_k^+/\delta_k^-)]+2\Bigr)+\log 3.
\\
&&
\nonumber
\check t := t+\sum_{k=1,4} 
\log\Bigl([\log_2(\delta_k^+/\delta_k^-)]+2\Bigr)+\log 3.
\end{eqnarray}
Then, with probability at least $1-e^{-t},$ for all $\delta_k\in [\delta_k^{-},\delta_k^{+}], k=1,2,3$
\begin{equation}
\label{first_bd}
\alpha_n(\delta_1,\delta_2)\leq 
2C_1 L(a)
{\mathbb E}\|\Xi\| (\sqrt{{\rm rank}(S)}\beta(S)\delta_1+\delta_2)+
4L(a)\delta_1\sqrt{\frac{\bar t}{n}}
+4L(a)a\frac{\bar t}{n},
\end{equation}
\begin{equation}
\tilde \alpha_n(\delta_1,\delta_2,\delta_3)\leq 
2C_2 L(a)
{\mathbb E}\|\Xi\| (\delta_2+\delta_3)+
4L(a)\delta_1\sqrt{\frac{\tilde t}{n}}
+4L(a)a\frac{\tilde t}{n},
\end{equation}
and 
\begin{equation}
\check \alpha_n(\delta_1,\delta_4)\leq 
2C_2 L(a){\mathbb E}\|\Xi\| \delta_4+
4L(a)\delta_1\sqrt{\frac{\check t}{n}}
+4L(a)a\frac{\check t}{n},
\end{equation}
where $C_1,C_2>0$ are numerical constants.
\end{lemma}

{\bf Proof}. 
We will prove in detail only the first bound (\ref{first_bd}). 
Talagrand's concentration inequality (in Bousquet's form, see Koltchinskii (2011), p. 25) implies that, for 
all $\delta_1, \delta_2>0,$ 
with probability at least $1-e^{-t}$
$$
\alpha_n(\delta_1,\delta_2)
\leq 2{\mathbb E}\alpha_n(\delta_1,\delta_2)+
2L(a)\delta_1\sqrt{\frac{t}{n}}
+4L(a)a\frac{t}{n}, 
$$
where we also used the bounds 
$$
|(\ell'\bullet f_{A})(f_{ A}-f_{S})|\leq 2 L(a)a,\ \ 
P(\ell'\bullet f_{A})^2(f_{ A}-f_{S})^2 \leq L^2(a)\|f_{ A}-f_{S}\|_{L_2(\Pi)}^2
\leq L^2(a)\delta_1^2
$$
that hold under the assumptions on the loss. The next step is to use standard 
Rademacher symmetrization and contraction inequalities (see, e.g., Koltchinskii (2011), sections 2.1, 2.2) to get 
\begin{equation}
\label{shest}
{\mathbb E}\alpha_n(\delta_1,\delta_2)\leq 16 L(a) 
{\mathbb E}\sup\{|R_n(f_A-f_S)|: A\in {\cal A}(\delta_1,\delta_2)\},
\end{equation}
where $R_n(f):=\sum_{j=1}^n \eps_j f(X_j),$ 
$\{\eps_j\}$ being 
i.i.d. Rademacher random variables independent of $\{(X_j,Y_j)\}$
and where we also used a simple fact that the Lipschitz constant of the function 
$u\mapsto \ell'(f_S+u)u$ is upper bounded by $4L(a).$ We will bound the 
expected sup-norm of the Rademacher process in the right hand side of (\ref{shest}).
Observe that
$$
R_n(f_A-f_S)=\langle \Xi, A-S\rangle, \ \ 
\Xi:= n^{-1}\sum_{j=1}^n \eps_j X_j,
$$
which implies 
\begin{eqnarray}
\label{ugu}
&&
|R_n(f_A-f_S)|\leq 
|\langle {\cal P}_L\Xi, {\cal P}_L(A-S)|+
|\langle \Xi, {\cal P}_L^{\perp}(A-S)| 
\\
&&
\nonumber
\leq \|{\cal P}_L\Xi\|_2\|{\cal P}_L(A-S)\|_2 
+
\|\Xi\|\|{\cal P}_L^{\perp}A\|_1
\\
&&
\nonumber
\leq 2\sqrt{2{\rm rank}(S)}\beta(S)\|\Xi\|\|f_A-f_S\|_{L_2(\Pi)}+
\|\Xi\|\|{\cal P}_L^{\perp}A\|_1,
\end{eqnarray}
where we used the facts that $A-S\in {\cal K}({\mathbb D};L;5)$
and also that 
$$
{\rm rank}({\cal P}_L \Xi)\leq 2 {\rm rank}(S),\ \ \|{\cal P}_L \Xi\|_2\leq 
2\sqrt{{\rm rank}({\cal P}_L \Xi)}\|\Xi\|. 
$$
Therefore, 
\begin{equation}
\label{ogo}
{\mathbb E}\sup\{|R_n(f_A-f_S)|: A\in {\cal A}(\delta_1,\delta_2)\}\leq 
{\mathbb E}\|\Xi\| (2\sqrt{2{\rm rank}(S)}\beta(S)\delta_1+\delta_2).
\end{equation}
It follows that with some numerical constant $C_1>0$ and with probability
at least $1-e^{-t},$
\begin{equation}
\alpha_n(\delta_1,\delta_2)\leq C_1 L(a)
{\mathbb E}\|\Xi\| (\sqrt{{\rm rank}(S)}\beta(S)\delta_1+\delta_2)+
2L(a)\delta_1\sqrt{\frac{t}{n}}
+4L(a)a\frac{t}{n}.
\end{equation}
We will make this bound uniform in $\delta_k\in [\delta_k^{-},\delta_k^{+}].$
To this end,
let $\delta_k^j:=\delta_k^+ 2^{-j}, 
j=0,\dots, [\log_2(\delta_{k}^{+}/\delta_k^{-})]+1.$ 
%Define 
%$\bar t = t+\sum_{k=1}^2 \log\Bigl([\log_2(\delta_k^+/\delta_k^-)]+2\Bigr)+\log 2.$
By the union bound, with probability at least $1-\frac{1}{3}e^{-t},$
for all $j_k=0,\dots, [\log_2(\delta_{k}^{+}/\delta_k^{-})]+1, k=1,2,$
\begin{equation}
\alpha_n(\delta_1^{j_1},\delta_2^{j_2})\leq C_1 L(a)
{\mathbb E}\|\Xi\| (\sqrt{{\rm rank}(S)}\beta(S)\delta_1^{j_1}+\delta_2^{j_2})+
2L(a)\delta_1^{j_1}\sqrt{\frac{\bar t}{n}}
+4L(a)a\frac{\bar t}{n},
\end{equation}
which implies that, 
for all $\delta_k\in [\delta_k^{-},\delta_k^{+}], k=1,2,$
\begin{equation}
\alpha_n(\delta_1,\delta_2)\leq 2C_1 L(a)
{\mathbb E}\|\Xi\| (\sqrt{{\rm rank}(S)}\beta(S)\delta_1+\delta_2)+
4L(a)\delta_1\sqrt{\frac{\bar t}{n}}
+4L(a)a\frac{\bar t}{n}.
\end{equation}

The proof of the second and the third bounds is similar. 
For instance, in the case of the second bound, the only difference is that instead 
of (\ref{ugu}) we use
\begin{equation}
\label{ugu'}
|R_n(f_A-f_S)|\leq \|\Xi\|(\|{\cal P}_L(A-S)\|_1+\|{\cal P}_L^{\perp}(A-S)\|_1),
\end{equation}
which yields (instead of (\ref{ogo}))
\begin{equation}
\label{ogo'}
{\mathbb E}\sup\{|R_n(f_A-f_S)|: A\in \tilde{\cal A}(\delta_1,\delta_2,\delta_3)\}\leq 
{\mathbb E}\|\Xi\| (\delta_2+\delta_3).
\end{equation}

\qed 

Note that 
\begin{equation}
\label{ta-ta}
(P-P_n)(\ell'\bullet f_{\hat S})(f_{\hat S}-f_{S})\leq 
\tilde \alpha_n(\|f_{\hat S}-f_S\|_{L_2(\Pi)};
\|{\cal P}_L^{\perp}\hat S\|_1;\|{\cal P}_L(\hat S-S)\|_1),
\end{equation}
\begin{equation}
\label{ta-ta-1}
(P-P_n)(\ell'\bullet f_{\hat S})(f_{\hat S}-f_{S})\leq 
\check \alpha_n(\|f_{\hat S}-f_S\|_{L_2(\Pi)};
\|\hat S-S\|_1),
\end{equation}
and also, if $\hat S-S\in {\cal K}({\mathbb D};L;b),$ then 
\begin{equation}
\label{tu-tu}
(P-P_n)(\ell'\bullet f_{\hat S})(f_{\hat S}-f_{S})\leq 
\alpha_n(\|f_{\hat S}-f_S\|_{L_2(\Pi)};\|{\cal P}_L^{\perp}\hat S\|_1).
\end{equation}
Assume for a while that 
\begin{equation}
\label{cond_delta}
\|f_{\hat S}-f_S\|_{L_2(\Pi)}\in [\delta_1^{-},\delta_1^{+}], 
\|{\cal P}_L^{\perp}\hat S\|_1\in [\delta_2^{-},\delta_2^{+}], 
\|{\cal P}_L(\hat S-S)\|_1\in [\delta_3^{-},\delta_3^{+}]. 
\end{equation}
First, we substitute (\ref{ta-ta-1}) in bound (\ref{tri})
and use the upper bound on $\check \alpha_n$ of Lemma \ref{main_lemma}. 
Observe also that, since $\hat V\in \partial\|\hat S\|_1,$ 
\begin{equation}
\label{sign''}
\langle \hat V, S-\hat S\rangle \leq \|S\|_1-\|\hat S\|_1.
\end{equation}
Therefore, we get 
\begin{eqnarray}
\label{chetyre_1_11}
&&
P(\ell \bullet f_{\hat S})+\frac{1}{2}\tau (a)\|f_{\hat S}-f_S\|_{L_2(\Pi)}^2
\\
&&
\nonumber
\leq 
P(\ell \bullet f_S)+ \eps (\|S\|_1-\|\hat S\|_1)+
\check \alpha_n(\|f_{\hat S}-f_S\|_{L_2(\Pi)};
\|\hat S-S\|_1)
\\
&&
\nonumber
\leq 
P(\ell \bullet f_S)+ 
\eps (\|S\|_1-\|\hat S\|_1)+
2C_2 L(a)
{\mathbb E}\|\Xi\| \|\hat S-S\|_1
\\
&&
\nonumber
+4L(a)\|f_{\hat S}-f_S\|_{L_2(\Pi)}\sqrt{\frac{\check t}{n}}
+4L(a)a\frac{\check t}{n}.
\end{eqnarray}
Assume that the constant $D$ in the condition on $\eps$ satisfies 
$D\geq 8C_2.$ Then, we have  
\begin{equation}
\label{eps_co}
\eps \geq DL(a)\Delta n^{-1/2}\geq 8C_2 L(a){\mathbb E}\|\Xi\|.
\end{equation}
Using the bound 
$$
4L(a)\|f_{\hat S}-f_S\|_{L_2(\Pi)}\sqrt{\frac{\check t}{n}}\leq 
\frac{1}{4}\tau (a)\|f_{\hat S}-f_S\|_{L_2(\Pi)}^2+
\frac{8L^2(a)}{\tau(a)}\frac{\check t}{n},
$$
we get from (\ref{chetyre_1_11})
\begin{eqnarray}
\label{chetyre_2_22}
&&
P(\ell \bullet f_{\hat S})
\leq 
P(\ell \bullet f_S)+ 
\eps (\|S\|_1-\|\hat S\|_1)
\\
&&
\nonumber
+\eps \|\hat S-S\|_1+
\biggl(\frac{8L^2(a)}{\tau(a)}
+4L(a)a\biggr)\frac{\check t}{n}
\\
&&
\nonumber
\leq 
P(\ell \bullet f_S)+ 
2\eps \|S\|_1+
\biggl(\frac{8L^2(a)}{\tau(a)}
+4L(a)a\biggr)\frac{\check t}{n}
\end{eqnarray}

We will now 
substitute (\ref{ta-ta}) in bound (\ref{chetyre})
and use the upper bound on $\tilde \alpha_n$ of Lemma \ref{main_lemma}. 
We will also bound $\langle {\rm sign}(S), S-\hat S\rangle$ as follows:
\begin{equation}
\label{sign}
|\langle {\rm sign}(S), S-\hat S\rangle|=
|\langle {\rm sign}(S), {\cal P}_L(S-\hat S)\rangle|
\leq \|{\rm sign}(S)\|\|{\cal P}_L(\hat S-S)\|_1
\leq \|{\cal P}_L(\hat S-S)\|_1.
\end{equation}
We get
\begin{eqnarray}
\label{chetyre_1}
&&
P(\ell \bullet f_{\hat S})+\frac{1}{2}\tau (a)\|f_{\hat S}-f_S\|_{L_2(\Pi)}^2
+ \eps \|{\cal P}_L^{\perp}(\hat S-S)\|_1
\\
&&
\nonumber
\leq 
P(\ell \bullet f_S)+ \eps \|{\cal P}_L(\hat S-S)\|_1+
\tilde \alpha_n(\|f_{\hat S}-f_S\|_{L_2(\Pi)};
\|{\cal P}_L^{\perp}\hat S\|_1;\|{\cal P}_L(\hat S-S)\|_1)
\\
&&
\nonumber
\leq 
P(\ell \bullet f_S)+ 
\eps \|{\cal P}_L(\hat S-S)\|_1+
2C_2 L(a)
{\mathbb E}\|\Xi\| (\|{\cal P}_L^{\perp}\hat S\|_1+\|{\cal P}_L(\hat S-S)\|_1)
\\
&&
\nonumber
+4L(a)\|f_{\hat S}-f_S\|_{L_2(\Pi)}\sqrt{\frac{\tilde t}{n}}
+4L(a)a\frac{\tilde t}{n}.
\end{eqnarray}
We still assume that $D\geq 8C_2$ and, thus, (\ref{eps_co}) holds.   
Using the bound 
$$
4L(a)\|f_{\hat S}-f_S\|_{L_2(\Pi)}\sqrt{\frac{\tilde t}{n}}\leq 
\frac{1}{4}\tau (a)\|f_{\hat S}-f_S\|_{L_2(\Pi)}^2+
\frac{8L^2(a)}{\tau(a)}\frac{\tilde t}{n},
$$
we get from (\ref{chetyre_1})
\begin{eqnarray}
\label{chetyre_2}
&&
P(\ell \bullet f_{\hat S})+\frac{1}{4}\tau (a)\|f_{\hat S}-f_S\|_{L_2(\Pi)}^2
+ \eps \|{\cal P}_L^{\perp}(\hat S-S)\|_1
\\
&&
\nonumber
\leq 
P(\ell \bullet f_S)+ 
\eps \|{\cal P}_L(\hat S-S)\|_1+
\frac{\eps}{4}(\|{\cal P}_L^{\perp}\hat S\|_1+\|{\cal P}_L(\hat S-S)\|_1)
\\
&&
\nonumber
\biggl(\frac{8L^2(a)}{\tau(a)}
+4L(a)a\biggr)\frac{\tilde t}{n}.
\end{eqnarray}
If 
$$
\biggl(\frac{8L^2(a)}{\tau(a)}
+4L(a)a\biggr)\frac{\tilde t}{n}
\geq \eps \|{\cal P}_L(\hat S-S)\|_1+
\frac{\eps}{4}(\|{\cal P}_L^{\perp}\hat S\|_1+\|{\cal P}_L(\hat S-S)\|_1),
$$
we conclude that 
\begin{equation}
\label{gu-gu}
P(\ell \bullet f_{\hat S})\leq 
P(\ell \bullet f_S)+ 
\biggl(\frac{16 L^2(a)}{\tau(a)}
+8L(a)a\biggr)\frac{\tilde t}{n},
\end{equation}
which suffices to prove the bound of the theorem. Otherwise, we use the assumption that 
$P(\ell\bullet f_{\hat S})> P(\ell\bullet f_S)$ to get 
the following bound from (\ref{chetyre_2}):
\begin{equation}
\nonumber
\eps \|{\cal P}_L^{\perp}(\hat S-S)\|_1\leq 
2\eps \|{\cal P}_L(\hat S-S)\|_1+
\frac{\eps}{2}(\|{\cal P}_L^{\perp}(\hat S-S)\|_1+\|{\cal P}_L(\hat S-S)\|_1).
\end{equation}
This yields
\begin{equation}
\nonumber
\frac{1}{2}\eps \|{\cal P}_L^{\perp}(\hat S-S)\|_1\leq 
\frac{5}{2}\eps \|{\cal P}_L(\hat S-S)\|_1,
\end{equation}
and, hence, $\hat S-S\in {\cal K}({\mathbb D};L;5).$ This fact allows us 
to use the bound on $\alpha_n$ of Lemma \ref{main_lemma}. We can modify 
(\ref{sign}) as follows 
\begin{eqnarray}
\label{sign_1}
&&
|\langle {\rm sign}(S), S-\hat S\rangle|=
|\langle {\rm sign}(S), {\cal P}_L(S-\hat S)\rangle|
\\
&&
\nonumber
\leq \|{\rm sign}(S)\|_2\|{\cal P}_L(\hat S-S)\|_2
\leq \sqrt{{\rm rank}(S)}\beta(S)\|f_{\hat S}-f_S\|_{L_2(\Pi)},
\end{eqnarray}
and, instead of (\ref{chetyre_1}), we get
\begin{eqnarray}
\label{chetyre_1'}
&&
P(\ell \bullet f_{\hat S})+\frac{1}{2}\tau (a)\|f_{\hat S}-f_S\|_{L_2(\Pi)}^2
+ \eps \|{\cal P}_L^{\perp}\hat S\|_1
\\
&&
\nonumber
\leq 
P(\ell \bullet f_S)+ 
\eps \sqrt{{\rm rank}(S)}\beta(S)\|f_{\hat S}-f_S\|_{L_2(\Pi)}+
\\
&&
\nonumber
2C_1 L(a)
{\mathbb E}\|\Xi\| (\sqrt{{\rm rank}(S)}\beta(S)\|f_{\hat S}-f_S\|_{L_2(\Pi)}+\|{\cal P}_L^{\perp}\hat S\|_1)+
\\
&&
\nonumber
+4L(a)\|f_{\hat S}-f_S\|_{L_2(\Pi)}\sqrt{\frac{\bar t}{n}}
+4L(a)a\frac{\bar t}{n}.
\end{eqnarray}
If $D\geq 2C_1,$ we have $\eps \geq 2C_1 L(a){\mathbb E}\|\Xi\|,$ and (\ref{chetyre_1'}) implies that 
\begin{eqnarray}
\label{chetyre_1''}
&&
P(\ell \bullet f_{\hat S})+\frac{1}{2}\tau (a)\|f_{\hat S}-f_S\|_{L_2(\Pi)}^2
\\
&&
\nonumber
\leq 
P(\ell \bullet f_S)+ 
\frac{3}{2\tau(a)}\beta^2(S){\rm rank}(S)\eps^2+ \frac{1}{6}\tau(a)\|f_{\hat S}-f_S\|_{L_2(\Pi)}^2+
\\
&&
\nonumber
\frac{3}{2\tau(a)}\beta^2(S){\rm rank}(S)\eps^2+\frac{1}{6}\tau(a)\|f_{\hat S}-f_S\|_{L_2(\Pi)}^2+
\\
&&
\nonumber
+\frac{24 L^2(a)}{\tau(a)}\frac{\bar t}{n}+\frac{1}{6}\tau(a)\|f_{\hat S}-f_S\|_{L_2(\Pi)}^2
+4L(a)a\frac{\bar t}{n}.
\end{eqnarray}
Therefore, we have 
\begin{eqnarray}
\label{chetyre_1'''}
&&
P(\ell \bullet f_{\hat S})
\leq 
P(\ell \bullet f_S)+ 
\frac{3}{\tau(a)}\beta^2(S){\rm rank}(S)\eps^2
+\biggl(\frac{24 L^2(a)}{\tau(a)}
+4L(a)a\biggr)\frac{\bar t}{n}.
\end{eqnarray}
The bound of the theorem will follow from (\ref{chetyre_2_22}), (\ref{gu-gu})
and (\ref{chetyre_1'''}) (provided that conditions (\ref{cond_delta}) hold).

We have to choose the numbers $\delta_k^{-}, \delta_k^{+}, k=1,2,3,4$ and 
establish the bound of the theorem when conditions (\ref{cond_delta})
do not hold. First note that, by the definition of $\hat S,$
$$
P_n(\ell\bullet \hat S)+ \eps \|\hat S\|_1 
\leq P_n (\ell \bullet 0)\leq Q,
$$
implying that $\|\hat S\|_1\leq \frac{Q}{\eps}.$
Next note that 
$$
\|{\cal P}_L^{\perp}\hat S\|_1=\|P_{L^{\perp}}\hat SP_{L^{\perp}}\|_1\leq \|\hat S\|_1\leq  \frac{Q}{\eps}
$$
and 
$$
\|{\cal P}_L(\hat S-S)\|_1\leq 2\|\hat S-S\|_1\leq  \frac{2Q}{\eps}+2\|S\|_1.
$$
Obviously, we also have 
$$
\|\hat S-S\|_1\leq \frac{Q}{\eps}+\|S\|_1.
$$
Finally, we have 
$
\|f_{\hat S}-f_S\|_{L_2(\Pi)}\leq 2a 
$
(since $\hat S, S\in {\mathbb D}$ and $\|f_{\hat S}\|_{L_{\infty}}\leq a, \|f_S\|_{L_{\infty}}\leq a$). Due to these facts, we can take 
$$
\delta_1^{+}:=2a,\ \delta_2^{+}:=\frac{Q}{\eps},\ \delta_3^{+}:=\frac{2Q}{\eps}+2\|S\|_1,
\delta_4^{+}:=\frac{Q}{\eps}+\|S\|_1,
$$
and, with this choice, $\delta_k^{+}, k=1,2,3, 4$ are upper bounds on the corresponding norms in (\ref{cond_delta}). We will also choose 
$$
\delta_1^{-}:= \frac{a}{\sqrt{n}},\ 
\delta_2^{-}:=\frac{L(a)a}{n\eps} \wedge (\delta_2^{+}/2),\ 
\delta_3^{-}:=\frac{L(a)a}{n\eps}\wedge (\delta_3^{+}/2),\ 
\delta_4^{-}:=\frac{L(a)a}{n\eps}\wedge (\delta_4^{+}/2).
$$
It is not hard to see that 
$$
\bar t \vee \check t\vee \tilde t \leq t(S;\eps)
$$
for a proper choice of numerical constant $B$ in the definition 
of $t(S;\eps).$ When conditions (\ref{cond_delta}) do not hold (which means that at least one of the numbers $\delta_k^{-}, k=1,2,3,4$ is not a lower bound on the corresponding 
norm), we still can use the bounds 
\begin{equation}
\label{ta-ta'}
(P-P_n)(\ell'\bullet f_{\hat S})(f_{\hat S}-f_{S})\leq 
\tilde \alpha_n(\|f_{\hat S}-f_S\|_{L_2(\Pi)}\vee\delta_1^{-};
\|{\cal P}_L^{\perp}\hat S\|_1\vee \delta_2^{-};\|{\cal P}_L(\hat S-S)\|_1\vee \delta_3^{-})
\end{equation}
\begin{equation}
\label{ta-ta-1'}
(P-P_n)(\ell'\bullet f_{\hat S})(f_{\hat S}-f_{S})\leq 
\check \alpha_n(\|f_{\hat S}-f_S\|_{L_2(\Pi)}\vee\delta_1^{-};
\|\hat S-S\|_1\vee \delta_4^{-})
\end{equation}
instead of (\ref{ta-ta}), (\ref{ta-ta-1}) and,  
in the case when $\hat S-S\in {\cal K}({\mathbb D};L;5),$ we can use the bound  
\begin{equation}
\label{tu-tu'}
(P-P_n)(\ell'\bullet f_{\hat S})(f_{\hat S}-f_{S})\leq 
\alpha_n(\|f_{\hat S}-f_S\|_{L_2(\Pi)}\vee \delta_1^{-};\|{\cal P}_L^{\perp}\hat S\|_1\vee \delta_2^{-})
\end{equation}
instead of bound (\ref{tu-tu}). It is easy now to modify the proof of 
(\ref{chetyre_1_11})--(\ref{chetyre_1'''}) to show that in this case 
we still have 
\begin{eqnarray}
\label{chetyre_1_1}
&&
\nonumber
P(\ell \bullet f_{\hat S})
\leq 
P(\ell \bullet f_S)+ 
\Bigl(\frac{3}{\tau(a)}\beta^2(S){\rm rank}(S)\eps^2\bigwedge 2\eps \|S\|_1\Bigr)
\\
&&
\nonumber
+C\biggl(\frac{L^2(a)}{\tau(a)}
\bigvee L(a)a\biggr)\frac{t(S;\eps)}{n},
\end{eqnarray}
which holds with probability at least $1-e^{-t}$ and implies 
the bound of the theorem. 

\qed

}

\end{document}